\documentclass[11 pt, a4 paper,twoside]{article}
\usepackage{amssymb,amsmath,graphicx,enumerate}
\usepackage{mathrsfs,longtable}
 \DeclareMathAlphabet{\mathpzc}{OT1}{pzc}{m}{it}
\numberwithin{equation}{section}
\newtheorem{thm}{Theorem}[section]
\newtheorem{lem}{Lemma}[section]

\newtheorem{defn}{Definition}[section]

\begin{document}
\begin{center}\textbf{LYAPUNOV TYPE INEQUALITY FOR HYBRID FRACTIONAL DIFFERENTIAL EQUATION WITH PRABHAKAR DERIVATIVE}\end{center}
\begin{center}
\textbf{Deepak B. Pachpatte, Narayan G. Abuj and Amol D. Khandagale}\\
Department of Mathematics,
\\Dr. Babasaheb Ambedkar Marathwada University,\\ Aurangabad - 431004, (M.S.), INDIA.\\
E-mail : {pachpatte@gmail.com, abujng@gmail.com and kamoldsk@gmail.com}
\end{center}
\begin{abstract}
 In this paper Lyapunov type inequality is developed for hybrid fractional boundary value problem involving the prabhakar fractional derivative.
\end{abstract}
${\textbf{\emph{MSC}}}: 26D10; 34B09; 33E12; 34A08$ 
\section{Introduction}
\paragraph{} In the fractional calculus the various integral inequalities plays very important role in the study of qualitative and quantitative
properties of solution of differential and integral equations.
The well-known Lyapunov result \cite{AML} sates that if the boundary value problem
\begin{equation}\label{1}
\left\{ {\begin{array}{*{20}{l}}{y^{\prime\prime}(t)}+q(t)y(t)= 0,\quad a<t< b,\quad\\
y(a)=y(b)=0, \end{array}} \right.
\end{equation}
has a nontrivial solution, where q(t) is real and continuous function, then\\
\begin{equation}\label{2}
\int_{a}^{b}|q(u)|du > \frac{4}{b-a}.
\end{equation}
\paragraph{}The study of Lyapunov inequalities for the fractional differential equation depends on a fractional differential operator involved and
it was initiated by Ferriera \cite{FRAC1}, also he derived a Lyapunov-type inequality for Riemann-Liouville fractional boundary value problem\\
\begin{equation}\label{3}
\left\{ {\begin{array}{*{20}{l}}{D^{\alpha}y(t)}+q(t)y(t)= 0,\quad a<t< b,\quad\\
y(a)=y(b)=0, \end{array}} \right.
\end{equation}
where $D^{\alpha}$ is the Riemann-Liouville fractional derivative of order $1<\alpha\leq2$ and $q(t):[a,b]\rightarrow \mathbb{R}$ is a continuous function. It has been proved that if $\eqref{3}$ has a nontrivial solution then\\
\begin{equation}\label{4}
\int_{a}^{b}|q(u)|du >\Gamma(\alpha)\bigg(\frac{4}{b-a}\bigg)^{\alpha-1}.
\end{equation}
 For $\alpha=2$ the inequality $\eqref{4}$ reduces to Lyapunov's classical inequality $\eqref{2}$.\\
 Also, Ferriera in  \cite{FRAC2} obtained a Lyapunov-type inequality for the Caputo fractional boundary value problem
\begin{equation}\label{5}
\left\{ {\begin{array}{*{20}{l}}{^{C}D^{\alpha}y(t)}+q(t)y(t)= 0,\quad a<t< b,\quad\\
y(a)=y(b)=0, \end{array}} \right.
\end{equation}
where ${^{C}D^{\alpha}}$ is the Caputo fractional derivative of order $1<\alpha\leq2$. It has been proved in \cite{FRAC2} that if $\eqref{5}$ has a nontrivial solution then
\begin{equation}\label{6}
\int_{a}^{b}|q(u)|du >\frac{\Gamma(\alpha){\alpha}^{\alpha}}{[(\alpha-1)(b-a)]^{\alpha-1}}.
\end{equation}
For $\alpha=2$, the inequality $\eqref{6}$ reduces to Lyapunov's classical inequality $\eqref{2}$.\\
  Also, Jleli and Samet \cite{MJ1,MJ2} modified the above inequalities for fractional differential equations with mixed boundary conditions.
\paragraph{} In \cite{SS}, Surang Sitho and et.al established Lyapunov type inequalities in two different cases for hybrid fractional boundary value problem
\begin{equation}\label{7}
\left\{ {\begin{array}{*{20}{l}}D_{a}^{\alpha}[\frac{y(t)}{f(t,y(t))}-\sum_{i=1}^{n}I_{a}^{\beta}h_{i}(t,y(t))]+q(t)y(t)= 0,\quad t\in (a,b), \\
y(a)=y^{\prime}(a)=y(b)=0. \end{array}} \right.
 \end{equation}
\,\,\, In $\eqref{7}$ $D_{a}^{\alpha}$ denotes the Riemann-Liouville fractional derivative of order $\alpha\in(2,3]$ starting from a point a, the function $y\in C([a,b],\mathbb{R}),g\in L^{1}((a,b],\mathbb{R}),$\\$f\in C^{1}([a,b]\times \mathbb{R},\mathbb{R}\setminus\{0\}),h_{i}\in C([a,b]\times \mathbb{R},\mathbb{R}),\forall i=1,2,...,n$ and $I_{a}^{\beta}$ is $\beta -th$ order Riemann-Liouville integral such that $\beta\geq\alpha $ with the lower limit at a point $a$.
\paragraph{} Recently, in \cite{SEAAA} the author's obtained the result on fractional diffential equation using the Prabhakar derivative
\begin{equation}\label{8}
  (\textbf{D}_{\rho,\mu,\omega,a+}^{\gamma}y)(t)+ q(t)y(t)= 0,\quad a<t<b,\quad 1<\mu\leq2,\quad \gamma,\rho,\omega \in {\mathbb{R^{+}}}
\end{equation}
with boundray conditions $ y(a)=y(b)=0$. Where $y \in C[a,b]$ and with the help Green function they obtained Lyapunov inequality for the fractional boundary value problem \eqref{8}
\paragraph{} Motivatied by  above work, in this paper we consider the following hybrid fractional differential equation involving the Prabhakar fractional derivative
\begin{equation}\label{9}
\left\{ {\begin{array}{*{20}{l}}\textbf{D}_{\rho,\mu,\omega,a+}^{\gamma}[\frac{y(t)}{p(t,y(t))}
-\sum_{i=1}^{n}\textbf{E}_{\rho,\mu,\omega,a+}^{\gamma}h_{i}(t,y(t))]+q(t)y(t)= 0,
\quad t\in(a,b), \\
y(a)=y(b)=0. \end{array}} \right.
\end{equation}
\,\,\,In $\eqref{9}$, $\textbf{D}_{\rho,\mu,\omega,a+}^{\gamma}$ denotes the Prabhakar derivative of order $\mu\in (1,2]$,\\
$ y\in C([a,b],\mathbb{R}), g\in L^{1}((a,b],\mathbb{R}), f\in C^{1}([a,b]\times \mathbb{R},\mathbb{R}\setminus\{0\}), h_{i}\in C([a,b]\times \mathbb{R},\mathbb{R}),\forall i=1,2,...,n$
and $\textbf{E}_{\rho,\mu,\omega,a+}^{\gamma}$ is the Prabhakar integral of order $\mu$ with lower limit at a point $a$. The Lyapunov type inequality is obtained for it.
\section{Preliminaries}
\begin{defn}
\cite{TRP} The generalized Mittag-Leffler function with three parameters is defined as,
\begin{equation}\label{a1}
E_{\rho,\mu}^{\gamma}(z)= \sum_{k=0}^\infty\frac{(\gamma)_kz^k}{\Gamma(\rho k+\mu)k!},\qquad \gamma, \rho, \mu \in\mathbb{C}, \Re(\rho)>0,
\end{equation}
where ${(\gamma)_k}$ is Pochhammer symbol defined by,\\
${(\gamma)_0}=1,\quad {(\gamma)_k}={\gamma(\gamma+1)...(\gamma+k-1)}$. for  $k=1,2,...$\\
For $\gamma=1$, the generalized Mittag-Leffler function \eqref{a1} reduces to the two-parameter Mittag-Leffler function given by
\begin{equation}\label{a2}
E_{\rho,\mu}(z):=E_{\rho,\mu}^{1}(z)=\sum_{k=0}^\infty\frac{z^k}{\Gamma(\rho k+\mu)},\qquad\rho, \mu \in \mathbb{C},\quad \Re(\rho)>0,
\end{equation}
and for $\mu=\gamma=1$, this function coincides with the classical Mittag-Leffler function $E_{\rho}(z)$
\begin{equation}\label{a3}
E_{\rho}(z):=E_{\rho, 1}^{1}(z)=\sum_{k=0}^\infty\frac{z^k}{\Gamma(\rho k+1)}, \qquad \rho\in \mathbb{C},\quad \Re(\rho)>0.
\end{equation}
Also, for $\gamma=0$ we have
$E_{\rho,\mu}(z)=\frac{1}{\Gamma(\mu)}$.
\end{defn}
\begin{defn}
\cite{RG} Let $f\in L^{1}[0,b]$, $0<x<b\leq\infty$, the prabhakar integral operator including generalized Mittag-Leffler function \eqref{a1} is defined as follows
\begin{equation}\label{a4}
\textbf{E}_{\rho,\mu,\omega,0+}^{\gamma}f(x)dx =\int_{0}^{x}(x-u)^{\mu-1}E_{\rho,\mu}^{\gamma}(\omega(x-u)^{\rho})f(u)du,\quad x>0
\end{equation}
where ${\rho,\mu,\omega,\gamma}\in \mathbb{C},$ with $\Re(\rho), \Re(\mu)>0.$\\
If for $\gamma=0$, the prabhakar integral operator coincides with the Riemann-Liouville fractional integral of order $\mu$;
$$\textbf{E}_{\rho,\mu,\omega,0+}^{0}f(x)=I_{0+}^{\mu}f(x),$$
where the Riemann-Lioville fractional integral is defined as
\begin{equation}\label{a5}
I_{0+}^{\mu}f(x)=\frac{1}{\Gamma(\mu)}\int_{0}^{x}(x-t)^{\mu-1}f(t)dt, \quad \mu\in \mathbb{C}, \Re(\mu)>0.
\end{equation}
\end{defn}
\begin{defn}
\cite{RG}  Let $f\in L^{1}[0,b]$, $0<x<b\leq\infty$, the Prabhakar derivative is defined as
\begin{equation}\label{a6}
\textbf{D}_{\rho,\mu,\omega,0+}^{\gamma}f(x)=\frac{d^m}{d{x}^m}\textbf{E}_{\rho, m-\mu,\omega,{0+}}^{-\gamma}f(x),
\end{equation}
where $\rho,\mu,\omega,\gamma\in \mathbb{C},$ with $\Re(\rho)>0$, $\Re(\mu)>o$, $m-1<\Re(\mu)<m$.\\
We note that the Prabhakar derivative genralizes the Riemann-Liouville fractional derivative
\begin{equation}\label{a7}
D_{0+}^{\mu}f(x)=\frac{d^m}{d{x}^m}\bigg(I_{0+}^{m-\mu}f \bigg)(x),\quad\mu \in\mathbb{C}, \Re(\mu)>0, m-1<\Re(\mu)<m.
\end{equation}
\end{defn}
\begin{lem}\label{L1} \cite{TRP} The Laplace transform of genralized Mittag-Leffler function $\eqref{a1}$ is given by
\begin{equation}\label{a8}
\mathscr{L}[x^{\mu-1}E_{\rho,\mu}^{\gamma}(\omega{x}^p)](s)=s^{-\mu}{(1-\omega s^{-\rho})}^{-\gamma},\qquad {|ws^{-\rho}|<1},
\end{equation}
for $\gamma,\rho,\mu,\omega, s\in \mathbb{C}, \Re(\mu)>0, \Re(s)>0$.
\end{lem}
\begin{lem}\label{L2} \cite{AAK} Let $\gamma, \rho, \mu, \omega, s \in \mathbb{C}$ with $\Re(\mu)>0.$ Then for any $ n \in \mathbb{N}$ differentiation of the genralized Mittag-Leffler function $\eqref{a1}$ is given by
\begin{equation}\label{a9}
\bigg(\frac{d}{dx}\bigg)^n[x^{\mu-1}E_{\rho,\mu}^{\gamma}(\omega x^\rho)]=x^{\mu-n-1}E_{\rho,\mu-n}^{\gamma}(\omega x^\rho).
\end{equation}
\end{lem}
\begin{lem}\label{L3}
\cite{SEAAA} The Laplace transform of Prabhakar integral $\eqref{a4}$ is given by
\begin{equation}\label{a10}
\mathscr{L}\{\textbf{E}_{\rho,\mu,\omega,0+}^{\gamma}f(x);s\}=s^{-\mu}(1-\omega s^{-\rho})^{-\gamma}F(s),
\end{equation}
where $F(s)$ is the Laplace transform of f(x), and it is written as
\begin{equation}\label{b1}
 F(s)= \mathscr{L}\{f(x);s\}=\int_{0}^{\infty}e^{-sx}f(x)dx,\,\quad s\in \mathbb{C}.
\end{equation}
\end{lem}
\begin{lem}\label{L4}
\cite{SEAAA} The laplace transform of Prabhakar derivative $\eqref{a6}$ is given by
\begin{equation}\label{b2}
\mathscr{L}\bigg\{\text{D}_{\rho,\mu,\omega,0+}^{\gamma}f(x);s\bigg\}= s^{\mu}(1-\omega s^{-\rho})F(s)-\sum_{k=1}^{m-1}s^{k}(\textbf{D}_{\rho,\mu-k-1,\omega,0+}^{\gamma}f)(0).
\end{equation}
\end{lem}
\begin{lem}\label{L5} \cite{SEAAA} If $f(x)\in C(a,b)\bigcap L(a,b),$ then
\begin{equation}\label{b3}
 \textbf{D}_{\rho,\mu,\omega,a+}^{\gamma}\textbf{E}_{\rho,\mu,\omega,a+}^{\gamma}f(x)=f(x),
\end{equation}
and if\\ $f(x),\textbf{D}_{\rho,\mu,\omega,a+}^{\gamma}f(x)\in C(a,b)\bigcap L(a,b)$ \text{then for} $c_{j} \in \mathbb{R},$ and $m-1<\mu\leq m$,\\
we have
\begin{align}\label{bb2}
\textbf{E}_{\rho,\mu,\omega,a+}^{\gamma} \textbf{D}_{\rho,\mu,\omega,a+}^{\gamma}f(x)&=f(x)+c_{1}(x-a)^{\mu-1}E_{\rho,\mu}^{\gamma}(\omega(x-a)^\rho){\nonumber}\\
&\quad +c_{2}(x-a)^{\mu-2}E_{\rho,\mu-1}^{\gamma}(\omega(x-a)^\rho)+...{\nonumber}\\
&\quad +c_{m}(x-a)^{\mu-m}E_{\rho,\mu-m+1}^{\gamma}(\omega(x-a)^\rho).
\end{align}
\end{lem}
\quad The authors had given the following proved lemma in \cite{SEAAA}.
\begin{lem}\label{T1}
 The Green function defined by
\begin{equation}\label{b4}
G(t,u)= \left\{
{\begin{array}{*{20}{l}}\frac{(t-a)^{\mu-1}E_{\rho,\mu}^{\gamma}(\omega(t-a)^\rho)}{(b-a)^{\mu-1}E_{\rho,\mu}^{\gamma}(\omega(b-a)^\rho)}{(b-u)^{\mu-1}E_{\rho,\mu}^{\gamma}(\omega(b-u)^\rho)}\\
 -(t-u)^{\mu-1}E_{\rho,\mu}^{\gamma}(\omega(t-u)^\rho),\quad\quad\quad\quad\quad\quad\quad\quad a\leq u\leq t\leq b, \\
\frac{(t-a)^{\mu-1}E_{\rho,\mu}^{\gamma}(\omega(t-a)^\rho)}{(b-a)^{\mu-1}E_{\rho,\mu}^{\gamma}(\omega(b-a)^\rho)}{(b-u)^{\mu-1}}E_{\rho,\mu}^{\gamma}(\omega(b-u)^\rho),\,\ a\leq t\leq u\leq b. \end{array}} \right.
\end{equation}
satisfies the following conditions :\\
1. For all $a\leq t, u\leq b, G(t,u)\geq 0.$\\
2. $\mathop {\max }\limits_{t\in[a,b]}G(t,u)= G(u,u)$, for $ u\in[a,b].$\\
3. The maximum of $G(u,u)$ is given at $u=\frac{a+b}{2}$ and has value
\begin{equation}\label{b5}
\mathop{\max}\limits_{u\in[a,b]}G(u,u)=G\bigg{(}\frac{a+b}{2},\frac{a+b}{2}\bigg{)}
= \bigg{(}\frac{b-a}{4}\bigg{)}^{\mu-1}\frac{E_{\rho,\mu}^{\gamma}(\omega(\frac{b-a}{2})^{\rho})
E_{\rho,\mu}^{\gamma}(\omega(\frac{b-a}{2})^{\rho})}{E_{\rho,\mu}^{\gamma}(\omega(b-a)^{\rho})}.
\end{equation}
\end{lem}
\section{Main Results}
\paragraph{} In this section, we have obtained Lyapunov type inequalities in two different cases:\\
(I)   $h_{i}{(t,y(t))}=0,\,\,i=1,2,...,n $ and \\
(II)  $h_{i}{(t,y(t))}\neq0,\,\, i=1,2,...,n.$ \\\\
\textbf{Case I :} $ h_{i}{(t,y(t))}=0,\,\,  i=1,2,...,n$\\
 Here, we consider the problem $\eqref{8}$ with $h_{i}(t,y(t))=0$, $\forall t\in[a,b]$, and for $\mu\in(1,2]$.\\
 \,\,\,We first construct a Green function for the following boundary value problem
\begin{equation}\label{b6}
 \left\{ {\begin{array}{*{20}{l}}\textbf{D}_{\rho,\mu,\omega,a+}^{\gamma}[\frac{y(t)}{p(t,y(t))}]+q(t)y(t)= 0,\quad 1<\mu\leq 2,\quad\gamma,\rho,\omega\in \mathbb{R}^{+},\\
y(a)=y(b)=0. \end{array}} \right.
\end{equation}
\begin{thm}\label{thm1}
Let $y\in AC([a,b],\mathbb{R})$ be a solution of $\eqref{b6}$. Then the function $y(t)$ satisfies the following integral equation
\begin{equation}\label{b7}
y(t)= p(t,y(t))\int_{a}^{b}G(t,u)q(u)y(u)du,
\end{equation}
where the Green function G(t,u) is given by $\eqref{b4}$.\\
\end{thm}
\textbf{Proof}: Operating $\textbf{E}_{\rho,\mu,\omega,a+}^{\gamma}$ on hybrid fractional differential equation $\eqref{b6}$ and using lemma $\eqref{L5}$ for real constant
$c_{1}$ and $\,c_{2}$ we have
\begin{align*}
\frac{y(t)}{p(t,y(t))}&= c_{1}(t-a)^{\mu-1}E_{\rho,\mu}^{\gamma}(\omega(t-a)^\rho)+ c_{2}(t-a)^{\mu-2}E_{\rho,\mu-1}^{\gamma}(\omega(t-a)^\rho)-\textbf{E}_{\rho,\mu,\omega,a+}^{\gamma}q(t)y(t)\\
&= c_{1}(t-a)^{\mu-1}E_{\rho,\mu}^{\gamma}(\omega(t-a)^\rho)+ c_{2}(t-a)^{\mu-2}E_{\rho,\mu-1}^{\gamma}(\omega(t-a)^\rho)\\
&\qquad-\int_{a}^{t}(t-u)^{\mu-1}E_{\rho,\mu-1}^{\gamma}(\omega(t-a)^\rho)q(y)y(u)du.
\end{align*}
Now, by employing the boundary conditions we obtain the value of $c_{1}$ and $c_{2}$ as follows
\begin{align*}
y(a)= 0
&\Leftrightarrow 0 = c_{2}(t-a)^{\mu-2}E_{\rho,\mu-1}^{\gamma}(\omega(t-a)^\rho)-\int_{a}^{a}
(a-u)^{\mu-1}E_{\rho,\mu}^{\gamma}(\omega(a-u)^\rho)q(u)y(u)du,\\
&\Leftrightarrow 0 = c_{2}(t-a)^{\mu-2}E_{\rho,\mu-1}^{\gamma}(\omega(t-a)^\rho),\\
&\Leftrightarrow c_{2}= 0,\\
\end{align*}
and \\
\begin{align*}
y(b)=0
&\Leftrightarrow 0 = c_{1}(b-a)^{\mu-1}E_{\rho,\mu}^{\gamma}(\omega(b-a)^\rho)- \int_{a}^{b}(b-u)^{\mu-1}E_{\rho,\mu}^{\gamma}(\omega(b-u)^\rho)q(u)y(u)du,\\
&\Leftrightarrow c_{1}={\frac{1}{(b-a)^{\mu-1}E_{\rho,\mu}^{\gamma}(\omega(b-a)^\rho)}} \int_{a}^{b}(b-u)^{\mu-1}E_{\rho,\mu}^{\gamma}(\omega(b-u)^\rho)q(u)y(u)du.
\end{align*}
Therefore the unique solution of $\eqref{b6}$ is written as follows
\begin{align*}
{\frac{y(t)}{p(t,y(t))}}&=\int_{a}^{t}\bigg{[}\frac{(t-a)^{\mu-1}E_{\rho,\mu}^{\gamma}(\omega(t-a)^\rho)}{(b-a)^{\mu-1}E_{\rho,\mu}^{\gamma}(\omega(b-a)^\rho)} {(b-u)^{\mu-1}E_{\rho,\mu}^{\gamma}(\omega(b-u)^{\rho})}\\
 &\quad -{(t-u)^{\mu-1}E_{\rho,\mu}^{\gamma}(\omega(t-u)^\rho)}\bigg{]} q(u)y(u)du\\
 &\quad +\int_{t}^{b}\bigg{[}\frac{(t-a)^{\mu-1}E_{\rho,\mu}^{\gamma}(\omega(t-a)^\rho)}{(b-a)^{\mu-1}
 E_{\rho,\mu}^{\gamma}(\omega(b-a)^\rho)} {(b-u)^{\mu-1}E_{\rho,\mu}^{\gamma}(\omega(b-u)^\rho)}\bigg{]}q(u)y(u)du,\\
y(t)&= p(t,y(t))\int_{a}^{b}G(t,u)q(u)y(u)du,
\end{align*}
where G(t,u)is given by $\eqref{b4}.$
\begin{thm}\label{thm2}
Let $\mathscr{B}= C[a,b]$ be the Banach space equipped with norm $\|y\|=\mathop{\sup}\limits_{t\in[a,b]}|y(t)|$ and nontrivial continuous solution of the hybrid fractional boundary value problem
\begin{equation*}
\left\{ {\begin{array}{*{20}{l}}\textbf{D}_{\rho,\mu,\omega,a+}^{\gamma}[\frac{y(t)}{p(t,y(t))}]+q(t)y(t)= 0,\qquad a<t<b, \\
y(a)=y(b)=0, \end{array}} \right.\\
\end{equation*}
exists, then
\begin{equation}\label{b8}
\frac{1}{\|p\|}{\bigg(\frac{4}{b-a}\bigg)^{\mu-1}}{\frac{E_{\rho,\mu}^{\gamma}(\omega(b-a)^\rho)}{E_{\rho,\mu}^{\gamma}
(\omega(\frac{b-a}{2})^\rho)E_{\rho,\mu}^{\gamma}(\omega(\frac{b-a}{2})^\rho)}}<\int_{a}^{b}|q(u)|du,
\end{equation}
\quad where q(t) is a real and continuous function.
\end{thm}
\textbf{Proof:} According to theorem \eqref{thm1}, a solution of the above fractional boundary value problem satisfies the integral equation
\begin{equation}
y(t)= p(t,y(t))\int_{a}^{b}G(t,u)q(u)y(u)du,
\end{equation}
which by applying the indicated norm on both sides of it,gives\\
\begin{align*}
&\|y\|\leq\|p\|\|y\|\mathop{\max}\limits_{t\in[a,b]}|G(t,u)|\int_{a}^{b}|q(u)|du,\\
&1\leq\|p\|\mathop{\max}\limits_{t\in[a,b]}|G(t,u)|\int_{a}^{b}|q(u)|du.\\
\end{align*}
Using the second property of the Green function in Lemma $\eqref{T1}$, we get desired inequality
\begin{equation*}
1<\|p\|\bigg{(}\frac{b-a}{4}\bigg{)}^{\mu-1}\frac{E_{\rho,\mu}^{\gamma}(\omega(\frac{b-a}{2})^\rho)E_{\rho,\mu}^{\gamma}(\omega(\frac{b-a}{2})^\rho)} {E_{\rho,\mu}^{\gamma}(\omega(b-a)^\rho)}\int_{a}^{b}|q(u)|du,\\
\end{equation*}
\begin{equation*}\label{b9}
\frac{1}{\|p\|}{\bigg(\frac{4}{b-a}\bigg)^{\mu-1}}{\frac{E_{\rho,\mu}^{\gamma}(\omega(b-a)^\rho)}{E_{\rho,\mu}^{\gamma}(\omega(\frac{b-a}{2})^\rho)
E_{\rho,\mu}^{\gamma}(\omega(\frac{b-a}{2})^\rho)}}<\int_{a}^{b}|q(u)|du. \quad\Box
\end{equation*}\\\\
Now, we consider the second case.\\
\textbf{Case II :} $h_{i}{(t,y(t))}\neq 0$,\quad $i=1,2,...,n$.\\
In this case, we construct Lyapunov type inequality for hybrid fractional boundary value problem $\eqref{9}$.
\begin{thm}\label{thm3}
Let $y\in AC[a,b]$ be a solution of $\eqref{9}$, then the function y(t) satisfy the following integral equation,\\
\begin{equation}\label{b10}
y(t)= p(t,y(t))\int_{a}^{b}G(t,y(t))\big[y(u)q(u)-\sum_{i=1}^{n}h_{i}(u,y(u))\big]du,
\end{equation}
where G(t,u) is Green function defined as in $\eqref{b4}$.
\end{thm}
\textbf{Proof:} Operating Prabhakar integral on $\eqref{8}$ we get\\
\begin{align*}
&\frac{y(t)}{p(t,y(t))}-\sum_{i=1}^{n}\textbf{E}_{\rho,\mu,\omega,a+}^{\gamma}h_{i}(t,y(t))+c_{1}(t-a)^{\mu-1}E_{\rho,\mu}^{\gamma}(\omega(t-a)^\rho)\\ &+c_{2}(t-a)^{\mu-2}E_{\rho,\mu-1}^{\gamma}(\omega(t-a)^\rho)+\textbf{E}_{\rho,\mu,\omega,a+}^{\gamma}q(t)y(t)=0.\\
\end{align*}
Rearranging the terms, we have \\
\begin{align}\label{b11}
\frac{y(t)}{p(t,y(t))}&=\sum_{i=1}^{n}\textbf{E}_{\rho,\mu,\omega,a+}^{\gamma}h_{i}(t,y(t))-\int_{a}^{t}(t-u)^{\mu-1}E_{\rho,\mu}^{\gamma}(\omega(t-u)^\rho)q(u)y(u)du\nonumber\\
&\qquad+c_{1}(t-a)^{\mu-1}E_{\rho,\mu}^{\gamma}(\omega(t-a)^\rho)+c_{2}(t-a)^{\mu-2}E_{\rho,\mu-1}^{\gamma}(\omega(t-a)^\rho).
\end{align}
Now, by employing the boundary conditions we can obtain the value of coefficients $c_{1}$ and $c_{2}$ as\\
\begin{align*}
y(a)= 0 \Leftrightarrow & \sum_{i=1}^{n}\textbf{E}_{\rho,\mu,\omega,a+}^{\gamma}h_{i}(a,y(a))= c_{2}(t-a)^{\mu-2}E_{\rho,\mu-1}^{\gamma}(\omega(t-a)^\rho),\\
\Leftrightarrow & c_{2}=\frac{\sum_{i=1}^{n}\int_{a}^{a}(t-u)^{\mu-1}\textbf{E}_{\rho,\mu,\omega,a+}^{\gamma}h_{i}(u,y(u))(\omega(t-u)^\rho)du}
{(t-a)^{\mu-2}}E_{\rho,\mu}^{\gamma}(\omega(t-a)^\rho)\\
\Leftrightarrow & c_{2}=0,
\end{align*}
and
\begin{align*}
y(b)=0 \Leftrightarrow \sum_{i=1}^{n}\bigg[\textbf{E}_{\rho,\mu,\omega,a+}^{\gamma}h_{i}(t,y(t))\bigg]_{t=b}
&=-\int_{a}^{b}(b-u)^{\mu-1}E_{\rho,\mu}^{\gamma}(\omega(b-u)^\rho)q(u)y(u)du\\&\quad +c_{1}(b-a)^{\mu-1}E_{\rho,\mu}^{\gamma}(\omega(b-a)^\rho),
\end{align*}
\begin{align*}\label{b12}
\qquad\quad \Leftrightarrow  c_{1}&=\frac{(b-a)^{\mu-1}}{E_{\rho,\mu}^{\gamma}(\omega(b-a)^\rho)}\bigg{\{}\int_{a}^{b}(b-u)^{\mu-1}
E_{\rho,\mu}^{\gamma}(\omega(b-u)^\rho)q(u)y(u)du\\\nonumber
&\quad -\sum_{i=1}^{n}\int_{a}^{b}(b-u)^{\mu-1}{E_{\rho,\mu}^{\gamma}(\omega(b-u)^\rho)}h_{i}(u,y(u))du\bigg{\}}.
\end{align*}
Substituting these value of $c_{1}$ and $c_{2}$ in equation $\eqref{b11}$ we get
\begin{align*}
 y(t)&=p(t,y(t))\bigg\{\sum_{i=1}^{n}\textbf{E}_{\rho,\mu,\omega,a+}^{\gamma}h_{i}(t,y(t))-\int_{a}^{t}(t-u)^{\mu-1}E_{\rho,\mu}^{\gamma}(\omega(t-u)^\rho)q(u)y(u)du\\
&\quad+\frac{(b-a)^{1-\mu}(t-a)^{\mu-1}}{E_{\rho,\mu}^{\gamma}(\omega(b-a)^\rho)}E_{\rho,\mu}^{\gamma}(\omega(t-a)^\rho)\bigg[\int_{a}^{b}(b-u)^{\mu-1}E_{\rho,\mu}^{\gamma}(\omega(b-u)^\rho)q(u)y(u)du\\
&\quad-\sum_{i=1}^{n}\int_{a}^{b}(b-u)^{\mu-1}{E_{\rho,\mu}^{\gamma}(\omega(b-u)^\rho)}h_{i}(u,y(u))du\bigg]\bigg\}\\
&= p(t,y(t))\bigg\{\sum_{i=1}^{n}\int_{a}^{t}(t-u)^{\mu-1}E_{\rho,\mu}^{\gamma}(\omega(t-u)^\rho)h_{i}(u,y(u))du\\
&\quad-\int_{a}^{t}(t-u)^{\mu-1}E_{\rho,\mu}^{\gamma}(\omega(t-u)^\rho)q(u)y(u)du\\
&\quad+\frac{(b-a)^{1-\mu}(t-a)^{\mu-1}}{E_{\rho,\mu}^{\gamma}(\omega(b-a)^\rho)}E_{\rho,\mu}^{\gamma}(\omega(t-a)^\rho)
 \bigg[\int_{a}^{b}(b-u)^{\mu-1}E_{\rho,\mu}^{\gamma}(\omega(b-u)^\rho)q(u)y(u)du\bigg]\\
&\quad-\frac{(b-a)^{1-\mu}(t-a)^{\mu-1}}{E_{\rho,\mu}^{\gamma}(\omega(b-a)^\rho)}E_{\rho,\mu}^{\gamma}(\omega(t-a)^\rho)
\bigg[\sum_{i=1}^{n}\int_{a}^{b}(b-u)^{\mu-1}E_{\rho,\mu}^{\gamma}(\omega(b-u)^\rho)h_{i}(u,y(u))du\bigg]\bigg\},\\
y(t)&=p(t,y(t))\bigg\{-\int_{a}^{t}(t-u)^{\mu-1}E_{\rho,\mu}^{\gamma}(\omega(t-u)^\rho)q(u)y(u)du\\
&\quad+\int_{a}^{t}\frac{(t-a)^{\mu-1}(b-u)^{\mu-1}E_{\rho,\mu}^{\gamma}(\omega(t-a)^\rho)E_{\rho,\mu}^{\gamma}(\omega(b-u)^\rho)} {(b-a)^{\mu-1}E_{\rho,\mu}^{\gamma}(\omega(b-a)^\rho)}q(u)y(u)du\\
&\quad+\int_{t}^{b}\frac{(t-a)^{\mu-1}(b-u)^{\mu-1}E_{\rho,\mu}^{\gamma}(\omega(t-a)^\rho)E_{\rho,\mu}^{\gamma}(\omega(b-u)^\rho)} {(b-a)^{\mu-1}E_{\rho,\mu}^{\gamma}(\omega(b-a)^\rho)}q(u)y(u)du\\
&\quad+\sum_{i=1}^{n}\int_{a}^{t}(t-u)^{\mu-1}E_{\rho,\mu}^{\gamma}(\omega(t-u)^\rho)h_{i}(u,y(u))du\\
&\quad-\frac{(b-a)^{1-\mu}(t-a)^{\mu-1}}{E_{\rho,\mu}^{\gamma}(\omega(b-a)^\rho)}E_{\rho,\mu}^{\gamma}(\omega(t-a)^\rho) \sum_{i=1}^{n}\int_{a}^{t}(b-u)^{\mu-1}E_{\rho,\mu}^{\gamma}(\omega(b-u)^\rho)h_{i}(u,y(u))du\\
&\quad-\frac{(b-a)^{1-\mu}(t-a)^{\mu-1}}{E_{\rho,\mu}^{\gamma}(\omega(b-a)^\rho)}E_{\rho,\mu}^{\gamma}(\omega(t-a)^\rho) \sum_{i=1}^{n}\int_{t}^{b}(b-u)^{\mu-1}E_{\rho,\mu}^{\gamma}(\omega(b-u)^\rho)h_{i}(u,y(u))du\bigg\},\\
y(t)&= p(t,y(t)\bigg[\int_{a}^{b}G(t,u)q(u)y(u)-\sum_{i=1}^{n}\int_{a}^{b}G(t,u)h_{i}(u,y(u))du\bigg],\\
y(t)&= p(t,y(t))\int_{a}^{b}G(t,u)\bigg[q(u)y(u)-h_{i}(u,y(u))\bigg]du,
\end{align*}
which is desired result.$\qquad\qquad\Box$ \\
\paragraph{} To prove our next result we use the following condition \\
$|q(u)y(u)-\sum_{i=1}^{n}h_{i}((u,y(u))|\leq K|q(u)\|y\|$.
\begin{thm}\label{thm4}
Let $\mathscr{B}= C[a,b]$ be the Banach space equipped with norm $ \|y\|= \mathop{\sup}\limits_{t\in[a,b]}|y(t)|$,\\
and a nontrivial continuous solution of the hybrid fractional boundary value problem $\eqref{8}$ exist, then
\begin{equation}\label{c2}
\frac{1}{K\|p\|}{\bigg(\frac{4}{b-a}\bigg)^{\mu-1}}{\frac{E_{\rho,\mu}^{\gamma}(\omega(b-a)^\rho)}{E_{\rho,\mu}^{\gamma}(\omega(\frac{b-a}{2})^\rho)
E_{\rho,\mu}^{\gamma}(\omega(\frac{b-a}{2})^\rho)}}< \int_{a}^{b}|q(u)|du,
\end{equation}
\end{thm}
where q(t) is real and continuous function.\\
\textbf{Proof}: In accordance with theorem \eqref{thm3}, a solution of the above hybrid fractional boundary value problem $\eqref{8}$, satisfies the integral equation \\
\begin{align*}
y(t)=p(t,y(t))\int_{a}^{b}G(t,u)\bigg[y(u)q(u)-\sum_{i=1}^{n}h_{i}(u,y(u))\bigg]du,\\
\end{align*}
which by applying the indicated norm on both sides of it, gives \\
\begin{align*}
\|y\|&\leq\|p\|\int_{a}^{b}|G(t,u)||y(u)q(u)-\sum_{i=1}^{n}h_{i}(u,y(u))|du,\\
\|y\|&\leq\|p\|\mathop{\max}\limits_{t\in[a,b]}|G(t,u)|\int_{a}^{b}|q(u)y(u)-\sum_{i=1}^{n}h_{i}(u,y(u))|du,\\
\|y\|&\leq K\|p\|\|y\|\mathop{\max}\limits_{t\in[a,b]}|G(t,u)|\int_{a}^{b}|q(u)|du.\\
\end{align*}
Using the second property of theorem $\eqref{T1}$, we get the desired inequality\\
\begin{align*}
\int_{a}^{b}|q(u)|du>\frac{1}{K\|p\|}{\bigg(\frac{4}{b-a}\bigg)^{\mu-1}}{\frac{E_{\rho,\mu}^{\gamma}(\omega(b-a)^\rho)}
{E_{\rho,\mu}^{\gamma}(\omega(\frac{b-a}{2})^\rho)
E_{\rho,\mu}^{\gamma}(\omega(\frac{b-a}{2})^\rho)}}.
\qquad\Box
\end{align*}\\\\\\\\

\end{document}